\documentclass{article}[11pt]
\usepackage{geometry}
\usepackage{enumerate}
\usepackage{latexsym,booktabs}
\usepackage{amsmath,amssymb, bm}
\usepackage{graphicx}
\usepackage{subfigure}

\usepackage{amsthm}
\usepackage{thmtools}
\usepackage{float}        
\usepackage{longtable}
\usepackage{multirow}
\usepackage{diagbox}
\usepackage{subfigure}
\usepackage{natbib}        
\usepackage{url}
\usepackage{textcomp}      
\usepackage{color}
\usepackage{subfiles}

\usepackage{lscape}

\usepackage{accents}
\newcommand{\dbtilde}[1]{\accentset{\approx}{#1}}

\renewcommand{\vec}[1]{\boldsymbol{#1}} 

\usepackage{siunitx} %
\usepackage[flushleft, online]{threeparttablex} %

\usepackage{algorithmicx}
\usepackage[ruled]{algorithm}
\usepackage{algpseudocode}
\algdef{SE}[DOWHILE]{Do}{doWhile}{\algorithmicdo}[1]{\algorithmicwhile\ #1}%

\usepackage{cases}

\renewcommand{\arraystretch}{1.5}

\theoremstyle{definition}

\renewcommand\thmcontinues[1]{\textbf{Continued}}

\usepackage{authblk}

\title{A hybrid inventory policy for non-stationary lot-sizing problem with lateral transshipment}

\author[a]{Xiyuan Ma\thanks{Corresponding author: Xiyuan.Ma@ed.ac.uk}}
\author[a]{Roberto Rossi}
\author[a]{Thomas Welsh Archibald}
\affil[a]{Business School, University of Edinburgh, Edinburgh, United Kingdom}

\date{}

\begin{document}
\maketitle

\begin{abstract}



This paper addresses the two-stocking locations single item non-stationary stochastic lot-sizing problem. The inventory level at each location is reviewed periodically. Items can be reordered and received from a common central warehouse and can also be transshipped laterally from the other location.

Lateral transshipment is assumed to be proactive to re-distribute the stock between two stocking locations. Therefore, the order of action in each period is: transshipping (if necessary), reordering (if necessary) and satisfying the demand at each location and each installation. The costs are imposed on transshipping, ordering, holding, and back-ordering. The key issue in such systems is to determine the quantity of the lateral transshipment between depots and the order quantities from the warehouse to both locations.

We formulate the problem via stochastic dynamic programming to minimise the expected total cost. Since the number of actions increases exponentially as the feasible quantities of transshipment and replenishment grow, we develop a two-stage dynamic programming to improve computation efficiency. 
A near-optimal policy against to this two-stage formulation is introduced based on a mixed integer linear programming and receding-horizon approach. numerical experiments are implemented to demonstrate the performance of two-stage model and the heuristic algorithm.\\

\end{abstract}

\section{Introduction}
A supply chain focuses on the core activities within organisations required to convert raw materials or component parts through to finished products or services. We look upstream to suppliers and the product flows among various retailers in a scale of inventory management, in which lot-sizing problems take a critical role by researching the timing and the quantity of replenishment \citep{silver1981operations}. In multi-location circumstances, the flow of item(s) takes place between warehouses and retailers as regular replenishment, as well as between retailers and retailers. The latter movements among stocking locations in the same echelon are called lateral transshipment.

\begin{figure}[H]
  \centering
  \includegraphics[width=0.9\textwidth]{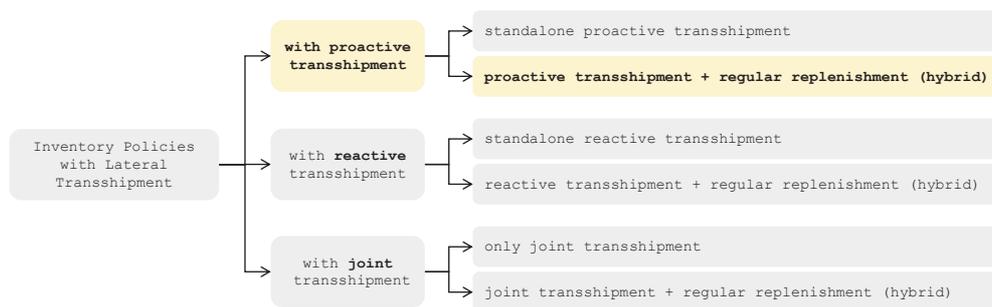}\\
  \caption{A classification of transshipment inventory policies.}\label{fig:classification}
\end{figure}

The inventory systems apply this lateral transshipment according to various aspects such as the products' natures, the structure of logistic routes, capital management policies, and so forth. If assume that the transshipment is always applied, the transshipment policy can be categorised as joint, hybrid and standalone with the consideration of the regular replenishment from the warehouse. \citet{paterson2011inventory} review the transshipment papers before 2010 and categorise transshipment into two types: proactive and reactive; the distinction is that the proactive transshipment is to reduce the risk of stockouts due to future demand while reactive transshipment is a recourse action to deal with existing stockouts or demand (also considered as an emergency replenishment). 
A classification is illustrated in Fig. \ref{fig:classification}, where we refer to the way that applies both proactive and reactive transshipment as a `joint' policy and that one applies both transshipment (in either type) and regular replenishment as a `hybrid' policy. Other criteria of classification are mentioned in \citep{paterson2011inventory}, for example, complete or partial pooling, centralised or decentralised systems.
Among vast transshipment literature, we focus on multi-location inventory systems with stochastic demand on a single item, where the system applies proactive transshipment plus the regular replenishment,as highlighted in Fig. \ref{fig:classification}.

The early research in this scope started from \citet{gross1963centralized}, showing that the optimal hybrid policy of a two-location single-period problem depends on the starting inventory level and cost parameters if the lead time is negligible, where the starting inventory levels divide the plane into six regions and each region has a corresponding optimal policy according to different cost parameters; the corresponding multi-location case is extended by \citet{karmarkar1977one} with a robust linear programming.
For a special case when the transshipment takes place at the end of the order cycle once only for non-identical locations, \citet{bertrand1998stock} show that it results in a single-period transshipment problem. The multi-period cases are tested to be of the similar characteristic as \citeauthor{gross1963centralized} in the optimal solutions by \citet{karmarkar1981multiperiod}.


\citet{hoadley1977two} state that applying proactive transshipment statically at the beginning of a period is beneficial to react to stock-outs as an additional opportunity is given.
In general, the transshipment timing has a significant influence on transshipment policy's performance, as \citet{tagaras2002effectiveness} explained, while it depends on system characteristics.
Most of the research on proactive transshipment is accomplished with or without regular replenishment by setting it at `static' points, referred to as \citet{agrawal2004dynamic}, where the re-distributions take place at the beginning of a period or a predefined point of a period.
A dynamic policy developed by \citet{agrawal2004dynamic} outperforms the static in terms of costs. They combine the real-time demand information to schedule transshipment with dynamic programming(DP) and solve near-optimal timing of transshipment and new stocking levels at stocking locations via a heuristic regarding the DP.

In replenishment models, the start or end (or any other point) of an order period provide `natural' opportunity to mitigate the mismatch between supply and demand by redistributing the stock over all locations. \citet{paterson2011inventory} argue that this explains the reason why the majority of research on proactive transshipment is based on a periodic review.

In this line of research, \citet{abouee2015optimal} provide a thorough analytical study and provide a `order-up-to curve' policy formed by four switching curves that divide the plane which is mathematically proved to be optimal, unlike the known optimal order-up-to level policy in multi-period inventory systems.

Considering other assumptions, \citet{diks1996controlling} consider the lead time and propose
the Consistent Appropriate Share Rationing (CAS) policy to balances the system stock to keep each location's fraction constant and the Balanced Stock (BS) rationing policy \citet{diks1998transshipments} as a general form of CAS, which is shown to be outperformed. They also showed that the system is benefited most when a large number of retailers, a high service level, or a long replenishment lead time is involved. 
The transshipment problem are also studied by \citet{archibald2010use}, \citet{powell2016perspectives} and \citet{meissner2018approximate} through approximate dynamic programming.

From the survey above, we notice that the influence of transship timing is emphasised by \citet{agrawal2004dynamic} for standalone proactive transshipment and by \citet{tagaras2002effectiveness} for hybrid policy, claiming that it is often set as a static point.
Meanwhile, the existing computation for hybrid policy in multi-period problems are either based on transshipment policies \citep{diks1996controlling,diks1998transshipments} that with strong assumption on system characteristics such as lead time, or initiated from systematical analysis \citep{abouee2015optimal} and apply approximate dynamic programming \citep{archibald2010use,meissner2018approximate}.
The gap lies in that no one has attempted to manipulate the formulation or to develop a solution approach from.
Besides, compared to stationary demand, there are relatively few studies in the literature that consider non-stationary demand, which, reflect the majority of practical circumstances.

In this paper, we propose a new optimisation method for 2-location stochastic lot-sizing problem with lateral transshipment under non-stationary demand. In our approach, we revise the feasible action space that integrates transshipment and replenishment. Inspired by \citet{agrawal2004dynamic}, we design the solution spaces as two splits from the original one according to the independence between transshipment and replenishment. Based on the separated solution spaces, we re-formulate the problem through a two-stage stochastic dynamic programming, where, for one period, the system first solves the order quantity and then the transshipping related decisions including the direction and the quantity.
To obtain the near-optimal solution that approximates the optimal expected total cost, we develop a heuristic solution approach in a receding-horizon framework, where we use an approximation technique that is different from \citeauthor{meissner2018approximate} for approximate dynamic programming. We make the following contributions to the stochastic lot-sizing literature.
\vspace*{-0.4\baselineskip}
\begin{itemize}\setlength{\itemsep}{-0.05cm}
\item We model the non-stationary stochastic lot-sizing problem with a hybrid policy combining proactive transshipment and regular replenishment by a two-stage stochastic dynamic formulation. The transshipping direction and quantity are determined in the first stage, and the order decisions and quantities are determined in the second stage as a single-location problem by considering the future cost that incorporates the whole system. 
\item We formulate the problem under a static uncertainty strategy in a mixed integer linear programming and develop a new heuristic algorithm to efficiently determine near-optimal policy parameters of the problem in the framework of receding horizon.
\item In a comprehensive numerical study, we show that two-stage formulation can well approximate the original dynamic programming and that the receding-horizon heuristic leads to tight optimality gaps.
\end{itemize}

The rest of this paper is structured as follows.
Section \ref{sec:description} describes the problem and formulates the problem as a stochastic dynamic programming.
Section \ref{sec:twoStageSDP} re-formulates the problem and proposes a two-stage formulation based on the one introduced in Section \ref{sec:generalSDP}. 
In Section \ref{sec:algorithm}, a heuristic algorithm is developed to compute near-optimal policy for the 2-location problems with lateral transshipment.
A computational analysis is presented in Section \ref{sec:computation}. And future researches are indicated in a discussion in Section \ref{sec:conclusion}.

\section{Problem Description}\label{sec:description}
We consider an inventory network of two stocking locations, at each of which demands occur randomly following the independent non-stationary stochastic distributions over a planning horizon of $T$ periods. The periods' demands $d_t^j$, for $t=1,\cdots,T$ and $j=1,2$, are independent random variables with known probability density functions $g_t^j(\cdot)$. We assume that each location has unlimited stocking capacity. Any demand that cannot be satisfied immediately is back-ordered.

The inventory level at each stocking location is reviewed periodically according to the periods of the non-stationary demand.
The supplement of stocks can be reordered and received from a common central warehouse and also be transshipped laterally from the other location in the network. We specify that only actual commodities (positive outstanding inventory level at one location) can be transshipped and so no back-order can be transferred.
For the replenishment, distinct locations across the system are of the same review periods and so locations are assumed to be replenished simultaneously if replenishment is applicable. And the warehouse is assumed to have sufficient inventory capacity to afford orders from stocking locations.
Regardless of the source, we assume that both supplements are issued and received instantaneously with negligible lead time.

At the beginning of a period of $t$, the inventory levels at two locations are both reviewed.
The system then makes decisions on transshipment and replenishment sequentially, including the transshipping direction, transshipping quantity and regular order quantity.
The cost incurred at every period over the planning horizon consists of
\begin{itemize}\setlength{\itemsep}{-0.1cm}
\item a transshipment cost $u(x) \triangleq R + vx$ for $x>0$, otherwise $u(x)=0$, where $x$ is the number of units transshipped, $R$ is the fixed cost per transshipment, and $v$ is the cost per unit transshipped;
\item an ordering cost $c(Q)\triangleq K + zQ$ for $Q>0$, otherwise $c(Q) = 0$, where $Q$ is the order quantity, $K$ is the fixed ordering cost per order from the warehouse, and $z$ is the cost per unit ordered \citep{scarf1959optimality};
\item a linear holding cost $h$ charged on every unit carried from one period to the next or the end of the planning horizon;
\item and a linear penalty cost $b$ charged on every unit back-ordered at the end of each period.
\end{itemize}

The objective of this study is to identify the optimal policy that integrates both transshipment and replenishment to minimise the expected total cost over a finite planning horizon and generate a well-performed approach for an arbitrary non-stationary stochastic demand series. We start from formulating the problem via a stochastic dynamic programme \citep{bellmandynamic} in the next section. Table \ref{sec:App.notations} can be found summarising all notations in Appendix \ref{tab:App.notations}.


\section{The stochastic dynamic programming formulation}\label{sec:generalSDP}
We model the problem as a stochastic dynamic programme and specify each component of the programme to present the final optimisation problem.

Under a periodic review, we apply a time period a time period $t\in\{1,\ldots,T\}$ for a $T-$period problem as the stage of the programming. The state of the system is observed at the beginning of a period and is described by two factors, $i^1$ and $i^2$, the inventory level in stocking location 1 and 2, respectively. We denote $\vec{i}\triangleq\langle i^1, i^2\rangle\in\mathcal{I}_t$, where $\mathcal{I}_t$ is the state space of stage $t$.

For each state in any stage, an action $\vec{a}_t\triangleq\langle W_t, Q_t^1, Q_t^2 \rangle\in\mathcal{A}_t$ indicates to schedule an lateral transshipment of $|W_t|$ units and then a replenishment order of $Q_t^1>0$ and $Q_t^2>0$ units for two locations, respectively, at the beginning of stage $t$ for the state $\vec{i}\in\mathcal{I}_t$. The transshipping quantity is dependent on the state $\vec{i}$, where $W_t\in\{W_t|\min\{0,-i^2\}\leq W_t\leq\max\{0,i^1\}\}$, a positive $W_t$ indicates a transshipment from location 1 to 2, and location 2 to 1 for the negative. A transshipment will only be deployed from location $j$ if $i^j>0$, for $j=1,2$. This indicates the transition probability between two states as
\begin{equation}
\text{Pr}(i^1 = r,i^2 = k|\vec{i},\vec{a}_t) = \left\{
\begin{array}{ll}
g_t^1(i^1+Q_t^1-r)\times g_t^2(i^2+Q_t^2-k), & \hbox{$i^1, i^2 \leq 0$,}\\
g_t^1(i^1 - W_t +Q_t^1-r)\times g_t^2(i^2 + W_t +Q_t^2-k), & \hbox{otherwise.} \\
\end{array}
\right.\notag
\end{equation}
where $\text{Pr}(i^1 = r,i^2 = k|\vec{i},\vec{a}_t)$ denotes the probability that state $\vec{i}\in\mathcal{I}_t$ transits to the state $\langle r,k\rangle$ under the action $\vec{a}_t\in\mathcal{A}_t$ in period $t$.

Therefore, the immediate cost consisting of holding and penalty costs, $f_t(\vec{i},\vec{a}_t)$, can be derived as
\begin{equation}
f_t(\vec{i},\vec{a}_t) \triangleq  \sum\nolimits_{j=1}^{2}\mathbb{E}[h\max(0,i^j\mp W_t + Q_t^j-d_t^j) + b\max(0,d_t^j - i^j\pm W_t - Q_t^j)],\label{eq:sdp.1.f}
\end{equation}
where the upper operation of `$\mp$' and `$\pm$' are applied for $j=1$.

Let $\vec{i}'$ represents the state after realising action $\vec{a}_t$ on state $\vec{i}_t\in\mathcal{I}_t$, where $i^{'1}=i^1-W_t+Q_t^1-d_t^1$ and $i^{'2}=i^2+W_t+Q_t^2-d_t^2$. The expected total cost over periods $t,\ldots,T$ with state $\vec{i}\in\mathcal{I}_t$ at the beginning of stage $t$ can be represented as $C_t(\vec{i})$, where
\begin{equation}
C_t(\vec{i}) =
\min_{\vec{a}_t\in\mathcal{A}_t}\mathop{} \{ u(|W_t|) + \sum\limits_{j=1}^{2}c(Q_t^j) + f_t(\vec{i},\vec{a}_t) + \mathbb{E}[C_{t+1}(\vec{i}')]\},\label{eq:sdp.1.finalC}
\end{equation}
and
\begin{equation}
C_{T+1}(\vec{i}) = 0 \label{eq:sdp.1.boundary}
\end{equation}
is the boundary condition. The optimisation problem therefore can be modelled as a stochastic dynamic program to solve $C_1(\vec{i})$ with $\vec{i}\in\mathcal{I}_1$. For convenience of notation, we denote this formulation as `SDP-1'. 

\section{A two-stage stochastic dynamic programming formulation}\label{sec:twoStageSDP}
Despite that an exact optimal policy could be provided by solving $C_1(\vec{i})$ according to Eq. \eqref{eq:sdp.1.finalC}, the exponential growth of action space causes redundancy in the computation and is also extremely time-consuming, not with mentioning a long planning horizon in the realistic setting.

This section exploits the formulation introduced in Section \ref{sec:generalSDP} and
%
%
develop a two-stage formulation to minimise the expected total cost, where the action space is decoupled by splitting transshipping actions and ordering actions.

The formulation continues using the notation of the stage $t$, state $\vec{i}$, action and action space $\vec{a}_t=\langle W_t, Q_t^1, Q_t^2 \rangle\in\mathcal{A}_t$
. The differences between the new and the former formulation are as follows.
\begin{itemize}\setlength{\itemsep}{-0.1cm}
\item We decouple the action space $\mathcal{A}_t = \mathcal{R}_t\times\mathcal{Q}_t$ for any period $t$, where $\mathcal{R}_t\subset\mathbb{Z}$ represents the space of feasible transshipping direction and quantity $W_t$, and $\mathcal{Q}_t\subset\mathbb{N}$ represents the space of feasible reorder quantity $\vec{q}_t\triangleq\langle Q_t^1, Q_t^2\rangle$ for depot 1 and 2, respectively.
\item We introduce a modified expected immediate cost $\tilde{f}_t$ of $\vec{q}_t$ and the state after transshipments as
\begin{equation}
\tilde{f}_t(\tilde{\vec{i}},\vec{q}_t) \triangleq  \sum\nolimits_{j=1}^{2}\mathbb{E}[h\max(0,\tilde{i}^j+Q_t^j-d_t^j) + b\max(0,d_t^j - \tilde{i}^j - Q_t^j)];\label{eq:sdp.2.f}
\end{equation}
where, for the clarity, $\tilde{\vec{i}}\in\widetilde{\mathcal{I}}_t$ denotes the state after transshipments and $\widetilde{\mathcal{I}}_t$ denotes the state space when the transshipment has completed.
\end{itemize}

Therefore, the expected total cost over periods $t,\ldots,T$ starting from state $\vec{i}\in\mathcal{I}_t$ can be described as
\begin{equation}
\widetilde{C}_t(\vec{i}) = \min_{W_t\in\mathcal{R}_t}\mathop{}\{u(|W_t|) + \dbtilde{C}(\langle i^1-W_t,i^2+W_t\rangle)\},\label{eq:sdp.2.transship}
\end{equation}
where
\begin{equation}
\dbtilde{C}(\tilde{\vec{i}}) = \min_{\vec{q}_t\in\mathcal{Q}_t}\mathop{}\{\sum_{j=1}^{2}c(Q_t^j) + \tilde{f}_t(\tilde{\vec{i}},\vec{q}_t) + \mathbb{E}[\widetilde{C}_{t+1}(\langle \tilde{i}^1+Q_t^1-d_t^1, \tilde{i}^2+Q_t^2-d_t^2\rangle)]\},\label{eq:sdp.2.end}
\end{equation}
and
\begin{equation}
\widetilde{C}_{T+1}(\vec{i}) =  0 \label{eq:sdp.2.boundary}
\end{equation}
is the boundary condition. The optimisation problem therefore can be modelled as a stochastic dynamic program to solve $\widetilde{C}_1(\vec{i})$ with $\vec{i}\in\mathcal{I}_1$.
Let $G_t(\vec{i})$ denote the expected cost over horizon ($t,T$) with no action taken in the first-leading period $t$. We denote this formulation as ``SDP-2''.

\section{An LP-based algorithm for proactive lateral transshipment problems under receding horizon}\label{sec:algorithm}
Although the optimal hybrid policy in this paper's scope can be obtained by enumerating all possible order and transshipping quantities as presented in Section \ref{sec:twoStageSDP}, the computation complexity increases exponentially as the planning horizon expands, and it becomes impractical to apply a dynamic programming. Instead of an exact solution, we compromise the optimality reasonably to obtain near-optimal solutions via various ways. This section introduces an algorithm to solve a near-optimal policy for this transshipment problems.

According to the problem description, the objective function minimises the expected total cost that comprises four types of costs, where 
ordering and transshipping costs are associated with decisions 
while holding and penalty costs are based on closing inventories dependent on the given stochastic demand distribution(s).

This brings the difficulty of modelling stochastic lot-sizing problems twofold.
The first challenge occurs in modelling the cost for the closing inventories of a period. Here we introduce the first-order loss function and its complementary, to capture penalty and holding costs, respectively, which play key roles in inventory control \citep{silver1998inventory}. Recent applications of loss functions in inventory control involve works by \citet{rossi2014piecewise}, \citet{rossi2015piecewise} and \citet{xiang2018computing}. And based on these works, loss function and its complementary can be approximated by two piecewise linear functions, and it leads to a linear programming (LP) for this problem.
The second challenge is to incorporate the demand uncertainty when building up the LP model. \citet{dural2020benefit} clearly indicate that both static and static-dynamic uncertainty feature very competitive optimality gaps and fully dominate the dynamic strategy under the receding horizon control \citep{bookbinder1988strategies}; this conclusion also applies to our problem. Therefore, we focus on featuring a simpler structure with static strategy, which results easier to implement in practice, rather than complexing an optimisation model, and introducing it to a receding horizon framework to solve and improve the solution by setting reasonable terminal conditions.

The development of the LP model and generic procedures of the algorithm are introduced in detail hereafter.

\subsection{LP model under a static uncertainty strategy}
Under a static uncertainty strategy, decisions on transship and order quantities over the planning horizon, $W_t$ and $Q_t^j$ for $t=1,\ldots,T$ and $j=1,2$, are all predefined according to demands and the closing inventories of the last period. We introduce $I_{t-1}^j$ to denote the closing inventory level at location $j$ in period $t-1$ and $W_t$ the transship quantity between location 1 and 2 at period $t$.

Consider a random variable $\omega$ and a scalar variable $x$. The first order loss function is defined as $\mathcal{L}(x,\omega) = \mathbb{E}[\max(\omega - x, 0)]$ and its complementary as $\widehat{\mathcal{L}}(x,\omega) = \mathbb{E}[\max(x - \omega, 0)]$. Let $\bar{H}_t^j$ and $\bar{B}_t^j$ represent the expected outstanding and back-ordering inventory at the end of period $t$ at location $j$ respectively, then
\begin{eqnarray}
\bar{H}_t^j &=& \widehat{L}(\bar{X}_t^j,\mathop{} d_{1\ldots t}^j),\notag\\
\bar{B}_t^j &=& L(\bar{X}_t^j,\mathop{} d_{1\ldots t}^j),\notag
\end{eqnarray}
where $d_{1\ldots t}^j$ denotes the convolution of random variable $d_k^j$, and
\begin{equation}
\bar{X}_t^j \triangleq \bar{I}_0^j + \sum_{k=1}^t(Q_k^j \mp W_k^j)\label{eq:cumulativeInventory}
\end{equation}
indicates the expected cumulative inventory available to satisfy demand up to period $t$ at location $j$. In this way, the expected total cost can be described as
\begin{equation}
\sum\limits_{t=1}^T\{u(|W_t|) + \sum_{j=1}^2[h\bar{H}_t^j + b\bar{B}_t^j + c(Q_t^j)] \} \notag
\end{equation}
to comprise ordering, transshipping, holding and penalty costs, where $u(|W_t|)$ and $c(Q_t^j)$ denote the transshipping cost and ordering cost at location $j$ in period $t$, respectively.

The constraints for the problem involve the flow balance equations and domains for decision variables. By taking expectations (denoted as $\bar{\cdot}$ ). The optimal policy under a static uncertainty strategy can be obtained by solving the following model, denoted as ``LP-1''.
\begin{alignat}{2}
\min \quad &  \sum\limits_{t=1}^T\{u(|W_t|) + \sum_{j=1}^2[h\bar{H}_t^j + b\bar{B}_t^j + c(Q_t^j)] \}, &\label{eq:MILP-obj}\\
\mbox{s.t.} \quad & \text{for }t=1\ldots,T\text{ and }j=1,2, &\quad&\notag\\
& \bar{H}_t^j = \widehat{L}(\bar{X}_t^j,\mathop{} d_{1\ldots t}^j),&\quad& \label{eq:MILP-H}\\
& \bar{B}_t^j = L(\bar{X}_t^j,\mathop{} d_{1\ldots t}^j),&\quad&\label{eq:MILP-B}\\
& \bar{I}_{t-1}^j \mp W_t + Q_t^j - \bar{d}_t^j = \bar{I}_t^j,&\quad&\label{eq:MILP-balance}\\
& W_t\leq\max\{0,\bar{I}_{t-1}^1\},&\quad& \label{eq:MILP-Wdomain1}\\
& W_t\geq\min\{0,-\bar{I}_{t-1}^2\},&\quad& \label{eq:MILP-Wdomain2}\\
& \bar{H}_t^j, \bar{B}_t^j \geq 0.&\quad&\label{eq:MILP-domain}
\end{alignat}
%
%
\citet{rossi2014piecewise} presents the piecewise linear upper and lower bounds for the first order loss function that can be immediately embedded in programming models. According to Lemma 3, 10 and 11 in \citep{rossi2014piecewise}, we generate the bounds for $\bar{H}_t^j$ and $\bar{B}_t^j$ under the static uncertainty strategy as
\begin{equation}
\bar{H}_t^j \geq \bar{X}_t^j\sum\limits_{k=1}^i p_k - \sum\limits_{k=1}^ip_k\mathbb{E}[d_{1\ldots t}|\Omega_k]\label{eq:piecewiseH}
\end{equation}
and
\begin{equation}
\bar{B}_t^j \geq \bar{X}_t^j\sum\limits_{k=1}^i p_k - \sum\limits_{k=1}^ip_k\mathbb{E}[d_{1\ldots t}|\Omega_k]-X_t^j + \bar{d}_{1\ldots t},\label{eq:piecewiseB}
\end{equation}
where $N$ disjoint adjacent subregions $\Omega_1, \Omega_2,\cdots, \Omega_N$ partition the domain of random variables $d_t^j$. These inequalities are included in the model to enable the computation for linear programming.

\subsection{A heuristic approach for transshipment problem based on static uncertainty strategy}
Model LP-1 as above does not stand alone in solving process but is embedded in a receding horizon framework. This subsection introduces a heuristic algorithm  for the two-location lot-sizing problem with lateral transshipment with the following procedures.

\vspace{-0.1cm}
\alglanguage{pseudocode}
\begin{algorithm}[H]\renewcommand\baselinestretch{1.2}\selectfont
\small
\caption{Computing near-optimal policy for non-stationary lot-sizing problem with lateral transshipment.}
\label{Algorithm}
\begin{algorithmic}[1]
\State \textbf{Input:} demand rates $d_t^j$;
cost parameters ($K$, $z$, $R$, $v$, $h$, $b$);
an opening state of inventory $\vec{i}_0$; a natural number $n=0$.
\State \textbf{Output:} a hybrid policy with transship and order quantities $\tilde{W}_t$ and $\tilde{Q}_t^j$ for $t=1,\ldots, T$ and $j=1,2$.

\Do
    \State Randomly generate two series $\hat{d}_t^j$ as demands at location 1 and 2 according to demand rates $d_t^j$ for $t=1,\ldots,T$;
    \For {$k = 1 \to T$}
        \State Solve model LP-1 and obtain solutions $\hat{W}_t$ and $\hat{Q}_k^j$ for $k=
        t,\ldots,T$ with opening inventory $\vec{i}_{k-1}$ given, and update $\tilde{Q}_k^j = \hat{Q}_k^j$;
        \State  Generate the transship quantity's domain $\mathcal{W}_k = \{W_k|\min\{0,-i_{k-1}^2\}\leq W_k\leq\max\{0,i_{k-1}^1\}\}$ with opening inventory $\vec{i}_{k-1}$, compute
        \begin{equation}
        \tilde{W}_k = \arg\min_{W_k\in\mathcal{W}_k}\dbtilde{C}(\vec{i}_{k-1})\label{eq:algorithm.W}
        \end{equation}
        by Eq.\eqref{eq:sdp.2.transship} -- \eqref{eq:sdp.2.end}, where $Q_k^j,\ldots,Q_T^j$, $W_{k+1},\ldots,W_T$ are substituted by $\hat{Q}_k^j,\ldots,\hat{Q}_T^j$, $\hat{W}_{k+1},\ldots,\hat{W}_T$ from LP-1, and $d_k^j,\ldots,d_T^j$ by $\hat{d}_k^j,\ldots,\hat{d}_T^j$;
        \State Compute inventory $\vec{i}_k$ with static demand $\hat{d}_k^j$ and decisions $\tilde{W}_k$ and $\tilde{Q}_k^j$;
    \EndFor
    \State Update the number of experiments $n=n+1$;
\doWhile
\State \hspace{0.5cm} The number of experiments $n$ is not sufficient large to contribute an $\alpha \times100\%$ confidence interval.
\end{algorithmic}
\end{algorithm}

We simulate the control policy (both transshipping and ordering) that obtained by implementing the aforementioned approach. The simulation iterates on the feasible initial states of inventory implement a stopping rule so as to achieve an estimation error of ±0.1\% of the expected total cost with 0.95 confidence probability (e.g. see \citep{dural2020benefit}).

\section{Computational experiments}\label{sec:computation}
This section presents a computational analysis to evaluate the accuracy of the proposed two-stage formulation in Section \ref{sec:twoStageSDP} and the performance of the MILP-based receding-horizon heuristic in Section \ref{sec:algorithm}. 
In Section \ref{sec:experiment1}, we consider a test set comprising small instances with 4 periods with Poisson demands and investigate the performance of the two-stage model against the optimal SDP in Section \ref{sec:generalSDP}.
In Section \ref{sec:experiment2}, we extend the planning horizon to 10 periods and alter the demands to Normal distributions; we compare the effectiveness of the heuristic against the two-stage SDP formulation. Note that in this experiment, two locations are not necessarily identical.

We refer to the percentage optimality gap of the expected total cost (ETC) as the measure to the comparison, which is computed according to $100\times (\text{ETC}_2 - \text{ETC}_1)/\text{ETC}_1$. We name three optimising approaches aforementioned by ``optimal-SDP'', ``two-stage-SDP'', and ``heuristic''. In the 4-period test set, $\text{ETC}_1$ and $\text{ETC}_2$ represents results by ``optimal-SDP'' and ``two-stage-SDP'', respectively. And in 10-period test set, they stand for results by ``two-stage-SDP'' and ``heuristic'', respectively.

For each test set, we consider three major parameters: fixed ordering cost $K$, fixed transshipping cost $R$, and penalty cost $b$. The other cost parameters are designed according to the scale of planning horizon lengths and demands and keep consistent within one set of experiments.
However, we predefine constraints $K>R$ to assume that the transshipping takes predecease over the ordering and $K\leq 2R$ to ensure the system would not order only once for the planning horizon without any other transshipment. Besides, we set $v<b$ to assume that the system would not leave unmet demand back-ordered even though the transshipment is reasonably worthwhile. These parameters will be elaborated in each subsection as below.

All computations are performed by a 4.0 ($1.90+2.11$) gigahertz Intel(R) Core(TM) i$7-8650$U CPU with 16.0 gigabytes of RAM in JAVA 1.8.0\_201.

\subsection{4-period test set with Poisson-distributed demands}\label{sec:experiment1}
This test set is designed to investigate the accuracy of two-stage formulation against the optimal SDP. We compare the difference between the expected total cost computed by the original SDP in Section \ref{sec:generalSDP} and by the two-stage SDP in Section \ref{sec:twoStageSDP} regarding to percentage optimality gaps.

\begin{figure}[H]
  \centering
  \includegraphics[width=0.7\textwidth]{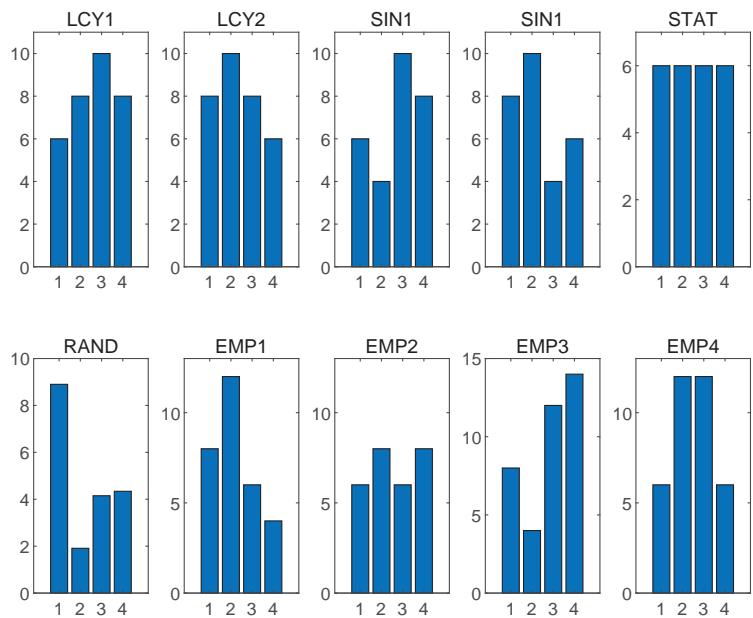}\\
  \caption{Demand patterns for 4-period instances.}\label{fig:testbed1}
\end{figure}

Due to the heavy computation for SDP approaches, we constrain that the demand means for Poisson distribution are no more than 20 as in Fig. \ref{fig:testbed1}. And we apply Latin hypercube sampling, which is introduced by \citet{mckay2000comparison}, to generate near-random sample of demand patterns from a multidimensional way.

We consider cost parameters including $K = 10, 20, 30$, $R = 5, 10, 20$ and $b = 3, 5$, which are implemented with two groups of unit ordering and transshipping pairs: $z=2$, $v=1$ and $z=1$, $v=0.5$ and $h=1$ for all experiments in this subsection. We also consider ten demand patterns: two life cycle patterns, one moves from the launch stage to maturity via a growth (LCY1) and the other moves from the growth stage through maturity and into decline (LCY2); two sinusoidal patterns, one with stronger (SIN1) and the other with weaker (SIN2) oscillations; a stationary demand pattern (STAT); a random demand pattern (RAND); and lastly, 4 empirical patterns derived according to \citet{strijbosch2011interaction}. Since the arbitrariness of locations, for one group of system parameters, two locations choose one demand out of ten sequentially without repetition of previous experiments, and demand patterns of two locations are not necessarily to be non-identical.

For the two-stage SDP approach, Table \ref{tab:4gap-parameter} reports the optimality gaps observed relative to the optimal SDP, pivoting demand patterns and various cost parameters. It should be noted that the two locations are of arbitrariness, therefore, we reformulate the results by pivoting the demand pattern of location 1 for the analysis.
The results for two-stage-SDP give the exact optimality gaps for these policies against the optimal SDP, which is on average 0.2924\%.
There is no obvious relation between optimality gaps and the variation in demand patterns or in cost parameters. The largest average arises under demand pattern LCY1 (0.3194\%), while it still does not deviate from the average.
We conclude that the two-stage SDP can nicely approximate the optimal SDP. We then proceed with the results by two-stage formulation as the benchmark in evaluating the heuristic approach's performance.

\begin{table}[H]\renewcommand\arraystretch{0.8}
\small
\centering
\caption{Average percent optimality gap over 4-period test set under pivoting parameters (\%).
}\label{tab:4gap-parameter}
\begin{tabular}{cp{1.5cm}<{\centering}}
\toprule
Pivoting parameters&  two-stage-SDP\\
\midrule
\textbf{demand patterns} &\\
LCY1&0.3194\\
LCY2&0.3076\\
SIN1&0.2930\\
SIN2&0.2863\\
STAT&0.2815\\
RAND&0.2886\\
EMP1&0.2890\\
EMP2&0.2850\\
EMP3&0.2896\\
EMP4&0.2840\\
\textbf{fixed ordering cost} \bm{$K$}&\\
10&0.2891\\
20&0.2836\\
30&0.3045\\
\textbf{fixed transshipping cost} \bm{$R$}&\\
5&0.2875\\
10&0.3054\\
20&0.2829\\
\textbf{penalty cost} \bm{$b$}&\\
3&0.2807\\
5&0.2829\\
\midrule
\textbf{Average}&0.2924\\
\bottomrule
\end{tabular}
\end{table}

\begin{figure}[H]
  \centering
  \includegraphics[width=0.8\textwidth]{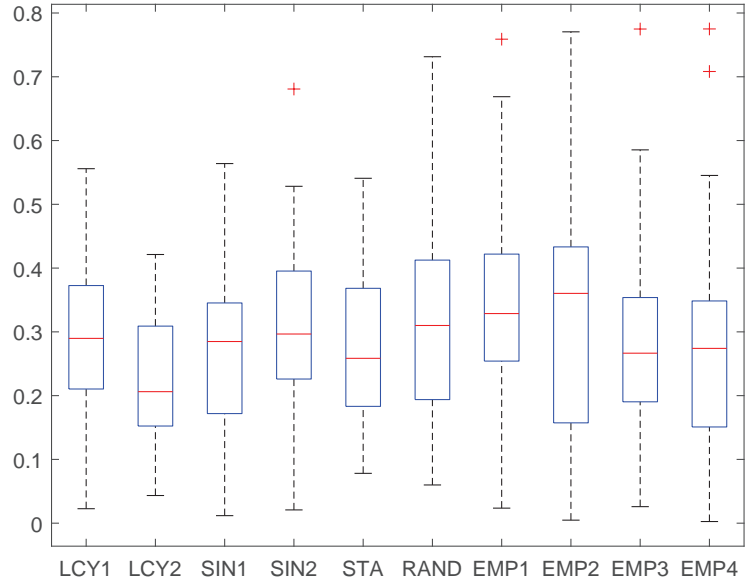}\\
  \caption{Box plots on the optimality gap regarding demand patterns.}\label{fig:test1_box_demand}
\end{figure}

\begin{figure}[H]
  \centering
  \includegraphics[width=0.8\textwidth]{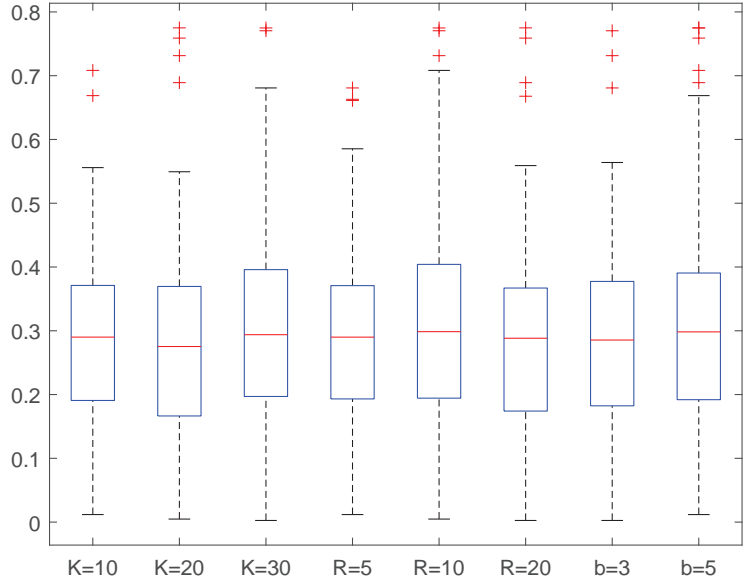}\\
  \caption{Box plots on the optimality gap regarding cost parameters.}\label{fig:test1_box_parameters}
\end{figure}

Focusing on characteristics of results inside, Fig. \ref{fig:test1_box_demand} and \ref{fig:test1_box_parameters} illustrate the distribution of optimality gaps against demand patterns and cost parameters, where all box plots present a Normal-distribution shape. More specifically, in the view of demand pattern, we observe that the mediums shift remarkably when changing demand patterns, while they can all be ranged in percentage of [0.2,0.4]. The long lower and upper whiskers indicate that the optimality gaps vary amongst each demand pattern.
In the view of cost parameters, we observe that box plots are of similar mediums at $0.30\%$ and distributed slightly different. More outliers are noticed comparing with Fig. \ref{fig:test1_box_demand} for demand patterns.
Considering an average optimality gap $0.2924\%$ from Table \ref{tab:4gap-parameter}, we conclude that the average accuracy of two-stage SDP can remain in a small value to facilitate further experiments.

\subsection{10-period test set with Normally-distributed demands}\label{sec:experiment2}
This subsection extends the planning horizon to 10 periods. The purpose of implementing this test set is evaluating the effectiveness of the proposed heuristic against the two-stage formulation. Since the computation of piecewise linearisation parameters consumes a large amount of computation time for long-horizon Poisson demand, this test set focuses on Normal demand, for which \citet{rossi2014piecewise} present pre-computed optimal partitioning coefficients. Note that we apply 10 partitions in piecewise approximation.

For demands in Normal distribution, we introduce the standard deviation $\sigma_t=\rho\cdot d_t$, where $d_t$ is the mean of a Normal distribution, and $\rho$ denotes the coefficient of variation of the demand, which remains fixed over time as prescribed by \citet{bollapragada1999simple}. We set $\rho=0.1$ in this set of experiments to ensure 1,000 of experiments can reflect unbiased results.

The Normal-distributed demand will be implemented on fixed ordering cost $K=10,20$, fixed transshipping cost $R=5,10$, penalty cost $b=3,5$, holding cost $h=1$, unit ordering cost $z=0.5$ and unit transshipping cost $v=1$, which construct 8 parameter groups. Due to the computation complexity, in this set, we downsize the number of demand patterns to only include life cycle (LCY), sinusoidal (SIN), stationary (STA), random (RAND), empirical (EMP), and they can still cover majority realistic demand situations. With arbitrary combinations of these patterns to each location, in total, one group of parameters will be applied on 25 pairs of demand patterns that two locations do not necessarily take the same, and 200 instances will be tested in this set.

\begin{figure}[H]
  \centering
  \includegraphics[width=0.9\textwidth]{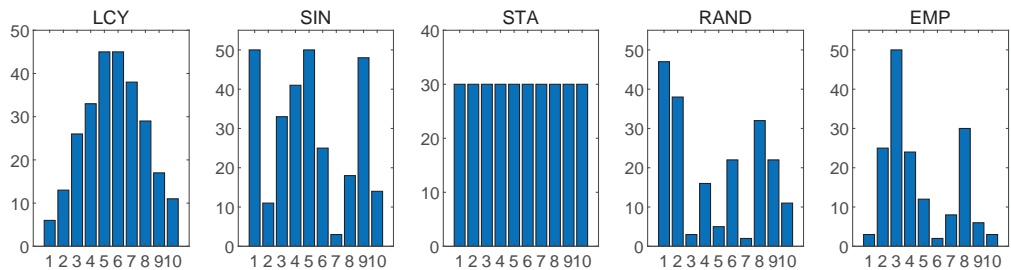}\\
  \caption{Demand patterns for 10-period instances.}\label{fig:testbed2}
\end{figure}

Table \ref{tab:10gap-parameter} reports the average percentage optimality gaps against expected total cost computed by two-stage SDP approach over 10-period instances pivoting five demand patterns and cost parameters. Similarly, we only pivot the demand pattern of location 1 for the analysis. 

This small set of instances reveal an average optimality gap $0.5515\%$. In the view of demand pattern, we observe that stationary pattern presents the minimum gap in $0.5140\%$, followed by random pattern ($0.5405\%$) and life cycle pattern ($0.5425\%$); sinusoidal pattern's average gap ($0.5557\%$) is slightly greater, and empirical pattern demonstrates the maximum ($0.6048\%$). Considering the shape of demand means in Fig. \ref{fig:testbed2}, we find that the performance of this receding-horizon heuristic deteriorates when demand is heavily non-stationary. 
In the view of cost parameters, the average gap pivoting fixed ordering cost $K$ varies more evidently than pivoting fixed transshipping cost $R$ and penalty cost $b$. And the gap increases when $K$ and $b$ independently raise but decreases when $R$ increases on the contrast.  

\begin{table}[H]\renewcommand\arraystretch{0.8}
\small
\centering
\caption{Average percent optimality gap over 10-period test set under pivoting parameters (\%).
}\label{tab:10gap-parameter}
\begin{tabular}{cp{5cm}<{\centering}}
\toprule
Pivoting parameters& MILP-based receding-horizon heuristic\\
\midrule
\textbf{demand patterns} &\\
LCY&0.5425\\
SIN&0.5557\\
STAT&0.5140\\
RAND&0.5405\\
EMP&0.6048\\
\textbf{fixed ordering cost} \bm{$K$}&\\
10&0.5373\\
20&0.5656\\
\textbf{fixed transshipping cost} \bm{$R$}&\\
5&0.5592\\
10&0.5438\\
\textbf{penalty cost} \bm{$b$}&\\
3&0.5506\\
5&0.5524\\
\midrule
\textbf{Average}&0.5515\\
\bottomrule
\end{tabular}
\end{table}

\section{Discussion}\label{sec:conclusion}
In this paper, we formulated and model a non-stationary 2-location, multi-period transshipment problem with back-order incases of stock-outs. We set the proactive transshipment takes place at the beginning of every period if applicable and established a stochastic dynamic programming to model the problem. However, the size of state and decision spaces makes it impossible to find the optimal policy for real-life sized problems, hence we introduced a two-stage stochastic dynamic programming, where the proactive transshipping quantity and order quantity are sequently determined.

Although the computation efficiency has been improved via the two-stage model, it is still hard to be applied. To obtain a near-optimal policy, we built a mixed integer linear programming under a static uncertainty strategy, which features a simpler structure and is increased by the receding-horizon framework, lead to a tight optimality gap shown in the experimental study.\\

Future research from this paper are originated from the current development or further extension on the approaches. They can be conducted on the following aspects.
\vspace*{-0.4\baselineskip}
\begin{itemize}\setlength{\itemsep}{-0.05cm}
\item The connection between stochastic dynamic programmings in Section \ref{sec:description} and \ref{sec:twoStageSDP}. The two-stage SDP is developed by treating transshipment and replenishment independent activities. From the numerical study in Section \ref{sec:experiment1}, we see that the difference between two formulations roughly follows a Normal distribution. At the same time, the average and the maxima can still remain in small values. A study on this relation can reveal where the difference originates.
\item The analytical study on the transshipment feasible space. In line 7 in the algorithm, to obtain the optimal transshipping quantity, one has to go through all feasible options, which is time-consuming even for a 10-period instance. An analytical study on the structure of the feasible space, which varies with the opening inventory at the current stage, could simplify the computation if a feature such as monotonicity, could be explored.
\item More generally, the study can be extended to involve more realistic assumptions; for example, (non-)identical lead times of replenishment and transshipment to different locations, lost-sale scheme to deal with unmet demand, and capacity imposed on either inventory storage or two means of transportation.
\item Other structures of the inventory system. This paper considers a two-echelon two-location problem transshipment problem. This simple structure can be extended to a more complex network to involve more echelons and connections, a system with $N$ warehouses and $M$ stocking locations, where $N$ and $M$ $\in\mathcal{N}$, in order to exploit the transshipment policy or condition in a general network.
\end{itemize}

\clearpage

\appendix

\linespread{1}

\renewcommand\thefigure{\Alph{section}\arabic{figure}}
\renewcommand\theequation{\Alph{section}\arabic{equation}}
\setcounter{figure}{0}
\setcounter{table}{0}
\setcounter{equation}{0}
\renewcommand\thetable{\Alph{section}\arabic{table}}

\section{Notations}\label{sec:App.notations}
\renewcommand\arraystretch{1.1}\selectfont\small
\begin{longtable}{rp{12cm}}
\caption{ Notations of important functions/parameters}\label{tab:App.notations}\\ \toprule
Functions& Explanation\\
\midrule
$d_j^t$ & demand for location $j$ in period $t$, $t=1,\ldots, T$;\\
$g_t^j(\cdot)$& probability density function of demand $d_j^t$;\\
$u(x)$ & transshipping cost, $u(x)=R+v(x)$ for $x>0$, $u(x)=0$ if $x\leq 0$;\\
$c(Q)$ & ordering cost, $c(Q) = K +zQ$ for $Q>0$, $c(Q)=0$ if $Q\leq 0$.\\
\midrule
$\vec{i}$ & a state in the stochastic dynamic programming that presents the inventory levels of two locations, $\vec{i}=\langle i^1, i^2\rangle\in\mathcal{I}_t$;\\
$\vec{a}_t$ & a feasible action for state $\vec{i}$ in the stochastic dynamic programming that integrates transshipping quantity and order quantities, $\vec{a}_t=\langle W_t, Q_t^1, Q_t^2\rangle\in\mathcal{A}_t$ and $W_t\in\{W_t|\min\{0,-i^2\}\leq W_t\leq\max\{0,i^1\}\}$;\\
$f_t(\vec{i},\vec{a}_t)$ & the expected holding and penalty cost of two locations that starts from state $\vec{i}$ applying action $\vec{a}_t$;\\
$C_t(\vec{i})$ & the expected total cost over periods $t,\ldots,T$ with state $\vec{i}\in\mathcal{I}_t$ at the beginning of stage $t$;\\
$G_t(\vec{i})$ & the expected cost over horizon ($t,T$) with no action taken in the first-leading period $t$.\\
\midrule
$\tilde{f}_t$ & the modified expected immediate cost of $\vec{q}_t$ and the state after transshipment;\\
$\widetilde{C}_t(\vec{i})$ & the expected total cost over periods $t$ to $T$ starting from state $\vec{i}\in\mathcal{I}_t$;\\
$\dbtilde{C}(\vec{i})$ & the minimised expected cost including the current ordering cost and immediate cost, and the future cost from period $t$ to $T$.\\
\bottomrule
\end{longtable}

\clearpage
\bibliographystyle{plainnat}
\bibliography{main}

\begin{thebibliography}{26}
\providecommand{\natexlab}[1]{#1}
\providecommand{\url}[1]{\texttt{#1}}
\expandafter\ifx\csname urlstyle\endcsname\relax
  \providecommand{\doi}[1]{doi: #1}\else
  \providecommand{\doi}{doi: \begingroup \urlstyle{rm}\Url}\fi

\bibitem[Abouee-Mehrizi et~al.(2015)Abouee-Mehrizi, Berman, and
  Sharma]{abouee2015optimal}
Hossein Abouee-Mehrizi, Oded Berman, and Shrutivandana Sharma.
\newblock Optimal joint replenishment and transshipment policies in a
  multi-period inventory system with lost sales.
\newblock \emph{Operations Research}, 63\penalty0 (2):\penalty0 342--350, 2015.

\bibitem[Agrawal et~al.(2004)Agrawal, Chao, and Seshadri]{agrawal2004dynamic}
Vipul Agrawal, Xiuli Chao, and Sridhar Seshadri.
\newblock Dynamic balancing of inventory in supply chains.
\newblock \emph{European Journal of Operational Research}, 159\penalty0
  (2):\penalty0 296--317, 2004.

\bibitem[Archibald et~al.(2010)Archibald, Black, and
  Glazebrook]{archibald2010use}
Thomas~W Archibald, Daniel~P Black, and Kevin~D Glazebrook.
\newblock The use of simple calibrations of individual locations in making
  transshipment decisions in a multi-location inventory network.
\newblock \emph{Journal of the Operational Research Society}, 61\penalty0
  (2):\penalty0 294--305, 2010.

\bibitem[Bellman(1957)]{bellmandynamic}
R~Bellman.
\newblock \emph{Dynamic Programming, Princeton, NJ, USA: Princeton University
  Press}.
\newblock 1957.

\bibitem[Bertrand and Bookbinder(1998)]{bertrand1998stock}
Louise~P Bertrand and James~H Bookbinder.
\newblock Stock redistribution in two-echelon logistics systems.
\newblock \emph{Journal of the Operational Research Society}, 49\penalty0
  (9):\penalty0 966--975, 1998.

\bibitem[Bollapragada and Morton(1999)]{bollapragada1999simple}
Srinivas Bollapragada and Thomas~E Morton.
\newblock A simple heuristic for computing nonstationary (s, s) policies.
\newblock \emph{Operations research}, 47\penalty0 (4):\penalty0 576--584, 1999.

\bibitem[Bookbinder and Tan(1988)]{bookbinder1988strategies}
James~H Bookbinder and Jin-Yan Tan.
\newblock Strategies for the probabilistic lot-sizing problem with
  service-level constraints.
\newblock \emph{Management Science}, 34\penalty0 (9):\penalty0 1096--1108,
  1988.

\bibitem[Diks and De~Kok(1998)]{diks1998transshipments}
EB~Diks and AG~De~Kok.
\newblock Transshipments in a divergent 2-echelon system.
\newblock In \emph{Advances in distribution logistics}, pages 423--447.
  Springer, 1998.

\bibitem[Diks and De~Kok(1996)]{diks1996controlling}
Erik~Bas Diks and AG~De~Kok.
\newblock Controlling a divergent 2-echelon network with transshipments using
  the consistent appropriate share rationing policy.
\newblock \emph{International Journal of Production Economics}, 45\penalty0
  (1-3):\penalty0 369--379, 1996.

\bibitem[Dural-Selcuk et~al.(2020)Dural-Selcuk, Rossi, Kilic, and
  Tarim]{dural2020benefit}
Gozdem Dural-Selcuk, Roberto Rossi, Onur~A Kilic, and S~Armagan Tarim.
\newblock The benefit of receding horizon control: Near-optimal policies for
  stochastic inventory control.
\newblock \emph{Omega}, 97:\penalty0 102091, 2020.

\bibitem[Gross(1963)]{gross1963centralized}
Donald Gross.
\newblock Centralized inventory control in multilocation supply systems.
\newblock \emph{Multistage inventory models and techniques}, 1:\penalty0 47,
  1963.

\bibitem[Hoadley and Heyman(1977)]{hoadley1977two}
Bruce Hoadley and Daniel~P Heyman.
\newblock A two-echelon inventory model with purchases, dispositions,
  shipments, returns and transshipments.
\newblock \emph{Naval Research Logistics Quarterly}, 24\penalty0 (1):\penalty0
  1--19, 1977.

\bibitem[Karmarkar(1981)]{karmarkar1981multiperiod}
Uday~S Karmarkar.
\newblock The multiperiod multilocation inventory problem.
\newblock \emph{Operations Research}, 29\penalty0 (2):\penalty0 215--228, 1981.

\bibitem[Karmarkar and Patel(1977)]{karmarkar1977one}
Uday~S Karmarkar and Nitin~R Patel.
\newblock The one-period, n-location distribution problem.
\newblock \emph{Naval Research Logistics Quarterly}, 24\penalty0 (4):\penalty0
  559--575, 1977.

\bibitem[McKay et~al.(1979)McKay, Beckman, and Conover]{mckay2000comparison}
Michael~D McKay, Richard~J Beckman, and William~J Conover.
\newblock A comparison of three methods for selecting values of input variables
  in the analysis of output from a computer code.
\newblock \emph{Technometrics}, 42\penalty0 (1):\penalty0 55--61, 1979.

\bibitem[Meissner and Senicheva(2018)]{meissner2018approximate}
Joern Meissner and Olga~V Senicheva.
\newblock Approximate dynamic programming for lateral transshipment problems in
  multi-location inventory systems.
\newblock \emph{European Journal of Operational Research}, 265\penalty0
  (1):\penalty0 49--64, 2018.

\bibitem[Paterson et~al.(2011)Paterson, Kiesm{\"u}ller, Teunter, and
  Glazebrook]{paterson2011inventory}
Colin Paterson, Gudrun Kiesm{\"u}ller, Ruud Teunter, and Kevin Glazebrook.
\newblock Inventory models with lateral transshipments: A review.
\newblock \emph{European Journal of Operational Research}, 210\penalty0
  (2):\penalty0 125--136, 2011.

\bibitem[Powell(2016)]{powell2016perspectives}
Warren~B Powell.
\newblock Perspectives of approximate dynamic programming.
\newblock \emph{Annals of Operations Research}, 241\penalty0 (1):\penalty0
  319--356, 2016.

\bibitem[Rossi et~al.(2014)Rossi, Tarim, Prestwich, and
  Hnich]{rossi2014piecewise}
Roberto Rossi, S~Armagan Tarim, Steven Prestwich, and Brahim Hnich.
\newblock Piecewise linear lower and upper bounds for the standard normal first
  order loss function.
\newblock \emph{Applied Mathematics and Computation}, 231:\penalty0 489--502,
  2014.

\bibitem[Rossi et~al.(2015)Rossi, Kilic, and Tarim]{rossi2015piecewise}
Roberto Rossi, Onur~A Kilic, and S~Armagan Tarim.
\newblock Piecewise linear approximations for the static--dynamic uncertainty
  strategy in stochastic lot-sizing.
\newblock \emph{Omega}, 50:\penalty0 126--140, 2015.

\bibitem[Scarf(1960)]{scarf1959optimality}
Herbert~E. Scarf.
\newblock {O}ptimality of (${s,S}$) policies in the dynamic inventory problem.
\newblock In K.~J. Arrow, S.~Karlin, and P.~Suppes, editors, \emph{Mathematical
  Methods in the Social Sciences}, pages 196--202. Stanford University Press,
  Stanford, CA, 1960.

\bibitem[Silver(1981)]{silver1981operations}
Edward~A Silver.
\newblock Operations research in inventory management: A review and critique.
\newblock \emph{Operations Research}, 29\penalty0 (4):\penalty0 628--645, 1981.

\bibitem[Silver et~al.(1998)Silver, Pyke, Peterson,
  et~al.]{silver1998inventory}
Edward~Allen Silver, David~F Pyke, Rein Peterson, et~al.
\newblock \emph{Inventory management and production planning and scheduling},
  volume~3.
\newblock Wiley New York, 1998.

\bibitem[Strijbosch et~al.(2011)Strijbosch, Syntetos, Boylan, and
  Janssen]{strijbosch2011interaction}
Leo~WG Strijbosch, Aris~A Syntetos, John~E Boylan, and Elleke Janssen.
\newblock On the interaction between forecasting and stock control: The case of
  non-stationary demand.
\newblock \emph{International Journal of Production Economics}, 133\penalty0
  (1):\penalty0 470--480, 2011.

\bibitem[Tagaras and Vlachos(2002)]{tagaras2002effectiveness}
George Tagaras and Dimitrios Vlachos.
\newblock Effectiveness of stock transshipment under various demand
  distributions and nonnegligible transshipment times.
\newblock \emph{Production and Operations Management}, 11\penalty0
  (2):\penalty0 183--198, 2002.

\bibitem[Xiang et~al.(2018)Xiang, Rossi, Martin-Barragan, and
  Tarim]{xiang2018computing}
Mengyuan Xiang, Roberto Rossi, Belen Martin-Barragan, and S~Armagan Tarim.
\newblock Computing non-stationary {($s, S$)} policies using mixed integer
  linear programming.
\newblock \emph{European Journal of Operational Research}, 271\penalty0
  (2):\penalty0 490--500, 2018.

\end{thebibliography}

\end{document}